\documentclass[10pt]{amsart}
            
\usepackage{amsmath}               
\usepackage{amsfonts}              
\usepackage{amsthm}                
\usepackage[utf8]{inputenc}
\usepackage{leftidx}
\usepackage{tikz-cd}
\usepackage{mathtools}

\newtheorem{tw}{Theorem}[section]
\newtheorem{lem}[tw]{Lemma}
\newtheorem{defi}[tw]{Definition}
\newtheorem{prop}[tw]{Proposition}
\newtheorem{rem}[tw]{Remark}
\newtheorem{cor}[tw]{Corollary}

\title{Cohomology of manifolds with structure group $U(n)\times O(s)$}
\author{Paweł Raźny}
\address{Institute of Mathematics \\
	Faculty of Mathematics and Computer Science \\
	Jagiellonian University in Cracow
	}
\email{pawel.razny@uj.edu.pl}


\begin{document}

\begin{abstract} We introduce a new spectral sequence for the study of $\mathcal{K}$-manifolds which arises by restricting the spectral sequence of a Riemannian foliation to forms invariant under the flows of $\{\xi_1,...,\xi_s\}$. We use this sequence to generalize a number of theorems from $K$-contact geometry to $\mathcal{K}$-manifolds. Most importantly we compute the cohomology ring and harmonic forms of $\mathcal{S}$-manifolds in terms of primitive basic cohomology and primitive basic harmonic forms (respectively). As an immediate consequence of this we get that the basic cohomology of $\mathcal{S}$-manifolds are a topological invariant. We also show that the basic Hodge numbers of $\mathcal{S}$-manifolds are invariant under deformations. Finally, we provide similar results for $\mathcal{C}$-manifolds. 
\end{abstract}
\maketitle
\section{Introduction}
In \cite{Yano,Yano2} a study of manifolds with a tensor field $f$ of type $(1,1)$ was initiated. The importance of the tensor field $f$ stems from the fact that its existence is equivalent to the reduction of the structure group of the manifold to $U(n)\times O(s)$. As such this generalizes both the concept of almost complex and almost contact manifolds. The properties of the curvature of such manifolds where further studied in \cite{B}. In particular, a special class of $f$-structures, called $\mathcal{S}$-structure, was introduced which generalizes K\"ahler and Sasakian manifolds. Since then these structures as well as their various generalizations (e.g. $\mathcal{K}$-structures, $K$-$f$-contact manifolds) were thoroughly studied.
\newline\indent The purpose of this article is to study the cohomological properties of such manifolds. In particular we are interested in the relations between the cohomology of the manifold and the basic cohomology of the foliation defined by the $\mathcal{S}$-structure. This study is motivated by the important role played by basic cohomology in Sasakian Geometry (e.g. the Sasakian version of the Calabi-Yau Theorem) as well as new results from \cite{Vais,GL}. We approach the problem by introducing and studying a new spectral sequence which is a variation of the spectral sequence of a Riemannian foliation (studied in e.g \cite{ss1,ss2,ss3}) and relates basic cohomology of a (almost) $\mathcal{K}$-manifold to its de Rham cohomology. Using this sequence we generalize a well known fact that for $K$-contact manifolds satisfying the hard Lefschetz property (equivalently the transverse hard Lefschetz property) the cohomology of the manifold in degree $r\leq n$ is isomorphic to the primitive basic cohomology (see \cite{lin}). Due to Poincar\'e duality this allows us to recreate the cohomology of the manifold from basic primitive cohomology and viceversa, which implies that in such a case basic cohomology is a topological invariant. In fact, our results applie to a more general class of manifolds, namely almost $\mathcal{S}$-structures satisfying the basic hard Lefschetz property. Using similar methods we also proof an analogous result for $\mathcal{C}$-structures which constitutes another special class of manifolds distinguished in \cite{B} and generalising the notion of quasi-Sasakian manifolds with a closed $1$-form $\eta$. An immediate corollary is that for manifolds with such structures the basic cohomologies are a topological invariant. In particular, this is true for any $\mathcal{S}$-manifolds and $\mathcal{C}$-manifolds. 
\newline\indent We provide two additional applications of the above results. Firstly, we classify Harmonic forms on $\mathcal{S}$-manifolds and $\mathcal{C}$-manifolds in terms of basic harmonic forms. This can be treated as a generalization of Proposition 7.4.13 from \cite{BG}. Secondly, we show that basic Hodge numbers of almost $\mathcal{S}$-manifolds and almost $\mathcal{C}$-manifolds which have the transverse hard Lefschetz property are invariant under deformations of such manifolds. This generalizes the main result from \cite{My2}. In particular, it is worth noting that this result applies to $K$-contact manifolds satisfying the hard Lefschetz property. Moreover, this partially answers Question 1.2 from \cite{My2}, in that it gives a positiva answer to that question for a new class of transversely  K\"ahler foliations.
\newline\indent We designate the subsequent section to preliminaries on basic cohomology and $\mathcal{S}$-structures. In Section \ref{Tor} we describe the aforementioned spectral sequence which will be the key tool in this article. We apply it in Sections \ref{MainS} and \ref{MainC} to prove the main results for almost $\mathcal{S}$-structures and almost $\mathcal{C}$-manifolds respectively. The final two sections contain the additional applications mentioned above.
\section{Preliminaries}
\subsection{Foliations}
We provide a quick review of transverse structures on foliations.
\begin{defi} A codimension q foliation $\mathcal{F}$ on a smooth n-manifold M is given by the following data:
\begin{itemize}
\item An open cover $\mathcal{U}:=\{U_i\}_{i\in I}$ of M.
\item A q-dimensional smooth manifold $T_0$.
\item For each $U_i\in\mathcal{U}$ a submersion $f_i: U_i\rightarrow T_0$ with connected fibers (these fibers are called plaques).
\item For all intersections $U_i\cap U_j\neq\emptyset$ a local diffeomorphism $\gamma_{ij}$ of $T_0$ such that $f_j=\gamma_{ij}\circ f_i$
\end{itemize}
The last condition ensures that plaques glue nicely to form a partition of M consisting of submanifolds of M of codimension q. This partition is called a foliation $\mathcal{F}$ of M and the elements of this partition are called leaves of $\mathcal{F}$.
\end{defi}
We call $T=\coprod\limits_{U_i\in\mathcal{U}}f_i(U_i)$ the transverse manifold of $\mathcal{F}$. The local diffeomorphisms $\gamma_{ij}$ generate a pseudogroup $\Gamma$ of transformations on T (called the holonomy pseudogroup). The space of leaves $M\slash\mathcal{F}$ of the foliation $\mathcal{F}$ can be identified with $T\slash\Gamma$.
\begin{defi}
 A smooth form $\omega$ on M is called basic if for any vector field X tangent to the leaves of $\mathcal{F}$ the following equality holds:
\begin{equation*}
i_X\omega=i_Xd\omega=0.
\end{equation*}
Basic 0-forms will be called basic functions henceforth.
\end{defi}
Basic forms are in one to one correspondence with $\Gamma$-invariant smooth forms on T. It is clear that $d\omega$ is basic for any basic form $\omega$. Hence, the set of basic forms of $\mathcal{F}$ (denoted $\Omega^{\bullet}(M\slash\mathcal{F})$) is a subcomplex of the de Rham complex of M. We define the basic cohomology of $\mathcal{F}$ to be the cohomology of this subcomplex and denote it by $H^{\bullet}(M\slash\mathcal{F})$. A transverse structure to $\mathcal{F}$ is a $\Gamma$-invariant structure on T. For example:
\begin{defi}
$\mathcal{F}$ is said to be transversely symplectic if T admits a $\Gamma$-invariant closed 2-form $\omega$ of maximal rank. $\omega$ is then called a transverse symplectic form. As we noted earlier $\omega$ corresponds to a closed basic form of rank q on M (also denoted $\omega$).
\end{defi}
\begin{defi}
$\mathcal{F}$ is said to be transversely holomorphic if T admits a complex structure that makes all the $\gamma_{ij}$ holomorphic. This is equivalent to the existence of an almost complex structure $J$ on the normal bundle $N\mathcal{F}:=TM\slash T\mathcal{F}$ (where $T\mathcal{F}$ is the bundle tangent to the leaves) satisfying:
\begin{itemize}
\item $L_XJ=0$ for any vector field $X$ tangent to the leaves.
\item if $Y_1$ and $Y_2$ are sections of the normal bundle then:
\begin{equation*}
 N_J(Y_1,Y_2):=[JY_1,JY_2]-J[Y_1,JY_2]-J[JY_1,Y_2]+J^2[Y_1,Y_2]=0
 \end{equation*}
where $[$ , $]$ is the bracket induced on the sections of the normal bundle (which can be defined by a choice of complement $N$ of $T\mathcal{F}$ via $\pi_N([$ , $])$).
\end{itemize}
\end{defi}
\begin{rem}
If $\mathcal{F}$ is transversely holomorphic we have the standard decomposition of the space of complex valued forms $\Omega^{\bullet}({M\slash\mathcal{F},\mathbb{C}})$ into forms of type (p,q) and $d$ decomposes into the sum of operators $\partial$ and $\bar{\partial}$ of order (1,0) and (0,1) respectively. Hence, one can define the Dolbeault double complex $(\Omega^{\bullet,\bullet}({M\slash\mathcal{F},\mathbb{C}}),\partial,\bar{\partial})$, the Fr\"{o}licher spectral sequence and the Dolbeault cohomology as in the manifold case. 
\end{rem}
\begin{defi}
$\mathcal{F}$ is said to be transversely orientable if T is orientable and all the $\gamma_{ij}$ are orientation preserving. This is equivalent to the orientability of $N\mathcal{F}$.
\end{defi}
\begin{defi}
$\mathcal{F}$ is said to be Riemannian if T has a $\Gamma$-invariant Riemannian metric. This is equivalent to the existence of a Riemannian metric g on $N\mathcal{F}$ with $L_Xg=0$ for all vector fields X tangent to the leaves.
\end{defi}
\begin{defi}
A foliation is said to be Hermitian if it is both transversely holomorphic and Riemannian.
\end{defi}
\begin{defi} A foliation $\mathcal{F}$ together with a triple $(g,J,\omega)$ consisting of a transverse Riemannian metric, transverse holomorphic structure and transverse symplectic form is called transversely K\"ahler if the following compatibility condition holds:
$$\omega(\cdot,\cdot)=g(J\cdot,\cdot)=\omega (J\cdot,J\cdot)$$ 
\end{defi}
We finish this section by recalling the spectral sequence of a Riemannian foliation.
\begin{defi} We put:
$$F^k_{\mathcal{F}}\Omega^r(M):=\{\alpha\in\Omega^r(M)\text{ }|\text{ } i_{X_{r-k+1}}...i_{X_1}\alpha=0,\text{ for } X_1,...,X_{r-k+1}\in\Gamma(T\mathcal{F})\}.$$
An element of $F^k_{\mathcal{F}}\Omega^r(M)$ is called an $r$-differential form of filtration $k$.
\end{defi}
The definition above in fact gives a filtration of the de Rham complex. Hence, via known theory from homological algebra we can cosntruct a spectral sequence as follows:
\begin{enumerate}
\item The $0$-th page is given by $E_0^{p,q}=F_{\mathcal{F}}^{p}\Omega^{p+q}(M)\slash F_{\mathcal{F}}^{p+1}\Omega^{p+q}(M)$ and $d^{p,q}_0:E_0^{p,q}\to E^{p,q+1}_0$ is simply the morphism induced by $d$.
\item The $r$-th page is given inductively by:
$$E_{r}^{p,q}:=Ker(d^{p,q}_{r-1})\slash Im(d^{p,q}_{r-1})=\frac{\{\alpha\in F^{p}_{\mathcal{F}}\Omega^{p+q}(M)\text{ }|\text{ } d\alpha\in F_{\mathcal{F}}^{p+r}\Omega^{p+q+1}(M)\}}{F^{p+1}_{\mathcal{F}}\Omega^{p+q}(M)+d(F_{\mathcal{F}}^{p-r+1}\Omega^{p+q-1}(M))}$$
\item The $r$-th coboundary operator $d_r:E_r^{p,q}\to E^{p+r,q-r+1}_r$ is again just the map induced by $d$ (due to the description of the $r$-th page this has the target specified above and is well defined).
\end{enumerate}
Furthermore, since the filtration is bounded this spectral sequence converges and its final page is isomorphic to the cohomology of the cochain complex (in this case the de Rham cohomology of $M$).
\begin{rem} The above spectral sequence can be thought of as a generalization of the Leray-Serre spectral sequence in de Rham cohomology to arbitrary Riemannian foliations (as opposed to fiber bundles).
\end{rem}
\subsection{Basic Hodge theory}
We devote this section to provide some background information on basic Hodge theory (see \cite{E1}) which will be applied in the final two section of this article. Firstly, we recall a special class of Riemannian foliations on which the aforementioned theory is greatly simplified:
\begin{defi} A codimension $q$ foliation $\mathcal{F}$ on a connected manifold $M$ is called homologically orientable if $H^{q}(M\slash\mathcal{F})=\mathbb{R}$. A foliation $\mathcal{F}$ on a manifold $M$ is called homologically orientable if its restriction to each connected component of $M$ is homologically orientable.
\end{defi}
We will later see in the subsequent section that all the foliation considered in this paper are in fact homologically orientable and hence we shall restrict our attention to this case throughout the rest of this subsection.
\newline\indent Let $\mathcal{F}$ be a homologically orientable Riemannian foliation on a manifold $M$. One can use the transverse Riemannian metric to define a basic Hodge star operator $*_b$ pointwise. This in turn allows us to define the basic adjoint operator:
$$\delta_b=(-1)^{q(r+1)+1}*_bd*_b.$$
\begin{rem} While we choose this to be the definition of $\delta_b$, it is in fact an adjoint of $d$ with respect to an appropriate inner product on forms induced by the transverse metric $g$. However, the definition of this inner product is quite involved and not necessary for our purpose. Although, we shall state some of the classical results of basic Hodge theory which use this inner product in their proof. For details see \cite{E1}.
\end{rem}
Using $\delta_b$ we can define the basic Laplace operator via:
$$\Delta_b=d\delta_b+\delta_bd.$$
As it turns out this operator has some nice properties similar to that of the classical Laplace operator. In particular, it is transversely elliptic in the following sense:
\begin{defi}
A basic differential operator of order $m$ is a linear map $D:\Omega^{\bullet}(M\slash\mathcal{F})\rightarrow\Omega^{\bullet}(M\slash\mathcal{F})$ such that in local coordinates $(x_1,...,x_p,y_1,...,y_q)$ (where $x_i$ are leaf-wise coordinates and $y_j$ are transverse ones) it has the form:
\begin{equation*}
D=\sum\limits_{|s|\leq m}a_s(y)\frac{\partial^{|s|}}{\partial^{s_1}y_1...\partial^{s_q}y_q}
\end{equation*}
where $a_s$ are matrices of appropriate size with basic functions as coefficients. A basic differential operator is called transversely elliptic if its principal symbol is an isomorphism at all points of $x\in M$ and all non-zero, transverse, cotangent vectors at x.
\end{defi}
In particular, this implies the following important result from \cite{E1}:
\begin{tw} Let $\mathcal{F}$ be a Riemannian homologically orientable foliation on a compact manifold $M$. Then:
\begin{enumerate}
\item $H^{\bullet}(M\slash\mathcal{F})$ is isomorphic to the space of basic harmonic forms $Ker(\Delta_b)$. In particular, it is finitely dimensional.
\item The basic Hodge star induces an isomorphism between $H^{k}(M\slash\mathcal{F})$ and $H^{q-k}(M\slash\mathcal{F})$ given by taking the class of the image through $*_b$ of a harmonic representative.
\end{enumerate}
\end{tw}
\subsection{Primitive basic cohomology}
Here we will recall some basic symplectic Hodge theory with main focus on primitive basic cohomology. To the best of our knowledge there is no concise source on the subject in its full generality (some special cases are treated in e.g. \cite{BG,lin}), hence for the readers convienience we provide a proof for the existence of basic primitive representatives and the Lefschetz decomposition on basic harmonic forms of transversely K\"ahler foliations.
\newline\indent Throughout this subsection let $M$ be a compact manifold endowed with a homologically orientable, transversely symplectic Riemannian foliation $\mathcal{F}$ of codimension $2n$. Firstly, let us note that by using the symplectic structure on the normal bundle $N\mathcal{F}$ we can define a symplectic star operator $*_s$ fiber by fiber in the standard way. This operator can in turn be used to define a number of other operators (on transverse forms i.e. sections of $\bigwedge^{\bullet} N^*\mathcal{F}$ which can be naturally identified with differential forms satisfying $i_{X}\alpha=0$ for $X\in\Gamma(T\mathcal{F})$) of interest:
$$L\alpha:=\omega\alpha,\quad \Lambda:=*_sL*_s,\quad d^{\Lambda}:=(-1)^{k+1}*_sd*_s=d\Lambda-\Lambda d.$$
\begin{rem} Note that if $(g,\omega,J)$ is a compatible triple consisting of a transverse Riemannian metric, transverse symplectic form and transverse almost complex structure, then by simple linear algebra:
$$J*_s=*_b,$$
where $J$ acts on a $k$-form via:
$$(J\alpha)(X_1,...,X_k)=\alpha(JX_1,...,JX_k).$$
this in turn allows one to compute that $\Lambda$ is dual to $L$ with respect to the metric $g$.
\end{rem}
We also denote by $L$ the morphism induced in cohomology by $L$. We recall the solution to the basic Bryliński conjecture from \cite{BC}:
\begin{tw}\label{BC} Let $M$ be a compact manifold endowed with a homologically orientable, transversely symplectic Riemannian foliation $\mathcal{F}$. Then the following conditions are equivalent:
\begin{enumerate}
\item (basic hard Lefschetz property) The map $L^k:H^{n-k}(M\slash\mathcal{F})\to H^{n+k}(M\slash\mathcal{F})$ is an isomorphism for all $k$.
\item Every basic cohomology class has a $d^{\Lambda}$-closed representative.
\end{enumerate}
\end{tw}
We move to some results on basic primitive forms required in this paper.
\begin{defi} A basic $(n-k)$-form $\alpha$ is said to be primitive if $L^{k+1}\alpha=0$ (or equivalently $\Lambda\alpha=0$) for $k\in\mathbb{N}$. Similarilly a basic cohomology class $[\alpha]$ of degree $(n-k)$ is said to be primitive if $L^{k+1}[\alpha]=0$. The space of basic primitive cohomology classes is denoted by $\mathcal{P}H^{\bullet}(M\slash\mathcal{F})$.
\end{defi}
The notion of primitive forms gives rise to the so called Lefschetz decomposition of forms:
\begin{prop}\label{prim0} Let $M$ be a compact manifold endowed with a homologically orientable, transversely symplectic Riemannian foliation $\mathcal{F}$. Let $\alpha$ be a basic $r$-form. Then $\alpha$ can be uniquely decomposed as:
$$\alpha=\sum\limits_{i} \omega^i\beta_{r-2i},$$
where $\beta_{r-2i}$ are basic primitive forms of degree $r-2i$ and are given by the formula:
$$\beta_{r-2i}=(\sum_k\frac{a_{i,k}}{k!}L^k\Lambda^{i+k}\alpha),$$
where $a_{i,k}$ are constants depending only on $(n,i,k)$.
\begin{proof} The proof of the unique decomposition and the explicit formula for $\beta_{r-2i}$ is well known from linear algebra (i.e. it is preciselly the same as in the manifold case). The fact that the forms $\beta_{r-2i}$ are basic follows from the explicit formula for $\beta_{r-2i}$.
\end{proof}
\end{prop}
The following two theorems show that if a manifold satisfies the basic hard Lefschetz property then the above decomposition descends to cohomology.
\begin{tw}\label{prim1} Let $M$ be a compact manifold endowed with a homologically orientable, transversely symplectic Riemannian foliation $\mathcal{F}$ satisfying the basic hard Lefschetz property. Let $\alpha$ be a basic cohomology class of degree $r$. Then $\alpha$ can be uniquely decomposed as:
$$\alpha=\sum\limits_{i} \omega^i\beta_{r-2i},$$
where $\beta_{r-2i}$ are basic primitive cohomology classes.
\begin{proof} Firstly, let us note that by the basic hard Lefschetz property for $i<n$ the mapping $L:H^i(M\slash\mathcal{F})\to H^{i+2}(M\slash\mathcal{F})$ is a monomorphism. In particular, this means that for $r\leq n$ we have:
$$H^r(M\slash\mathcal{F})=\mathcal{P}H^r(M\slash\mathcal{F})\oplus LH^{r-2}(M\slash\mathcal{F}).$$
Proceeding inductively we get:
$$H^r(M\slash\mathcal{F})=\bigoplus\limits_i L^i\mathcal{P}H^{r-2i}(M\slash\mathcal{F}).$$
which proves the theorem for $r\leq n$. For $r>n$ simply compose the above decomposition for $2n-r$ with $L^{r-n}$ and apply the hard Lefschetz property.
\end{proof}
\end{tw}
\begin{tw}\label{prim2} Let $M$ be a compact manifold endowed with a homologically orientable, transversely symplectic Riemannian foliation $\mathcal{F}$ satisfying the basic hard Lefschetz property. Every basic primitive cohomology class has a basic primitive representative.
\begin{proof} Let $\alpha$ be a $d^{\Lambda}$-closed representative of a given basic primitivie cohomology class. Then each primitive component of $\alpha$ is given by:
$$\beta_{r-2i}=(\sum_k\frac{a_{i,k}}{k!}L^k\Lambda^{i+k}\alpha),$$
as described earlier. By applying $d$ to the left hand side and noting that:
\begin{enumerate}
\item $d$ commutes with $L$,
\item $d$ commutes with $\Lambda$ up to $d^{\Lambda}$,
\item $d^{\Lambda}$ commutes with $\Lambda$,
\end{enumerate}
we see that each $\beta_{r-2i}$ is closed. We note that each component aside from $\beta_r$ has to be also exact as otherwise $[\alpha]$ would not be primitive. Hence, $[\beta_r]=[\alpha]$ which ends the proof.
\end{proof}
\end{tw}
Finally, we give a similar decomposition theorem for basic harmonic forms on transversely K\"ahler foliations.
\begin{tw}\label{prim3} Let $M$ be a compact manifold endowed with a homologically orientable, transversely K\"ahler foliation $\mathcal{F}$.  Let $\alpha$ be a basic harmonic $r$-form. Then $\alpha$ can be uniquely decomposed as:
$$\alpha=\sum\limits_{i} \omega^i\beta_{r-2i},$$
where $\beta_{r-2i}$ are basic harmonic forms which are primitive.
\begin{proof} We start by prooving that if $\alpha$ is a basic harmonic $k$-form then so is $L\alpha$. It is known (cf. \cite{E1}) that if a basic form on a transversely K\"ahler foliation is basic harmonic then it is in the kernel of the operators $\partial$, $\bar{\partial}$ and their adjoints (with respect to the transverse metric). In particular, it is in the kernel of the adjoint $(d^c)^*$ of the operator:
$$d^c:= i(\bar{\partial}-\partial)=J^{-1}dJ.$$
but we can compute similarilly as in the classical case:
$$d^{\Lambda}\alpha=(-1)^{k+1}*_sd*_s=(-1)^{k+1}*J^{-1}d*J^{-1}=(-1)^{k+1}*d^c*J^{-2}=(d^c)^{*}.$$
Hence, we have proved that $0=d^{\Lambda}\alpha=d\Lambda\alpha$ this together with $\delta\Lambda\alpha=\Lambda\delta\alpha=0$ implies that $\Lambda\alpha$ is harmonic. By adjointess the fact that $\Lambda$ preserves being harmonic implies that $L$ does so as well.
\newline\indent Now if $\alpha$ is a basic harmonic form representing a primitive class then the form $L^{n-k+1}\alpha$ is harmonic as well. However, since $L^{n-k+1}\alpha$ represents the trivial cohomology class it has to be equall to $0$. Hence, $\alpha$ is a primitive form itself. Now the theorem follows from Theorem \ref{prim1} and the conclusion of the previous paragraph.
\end{proof}
\end{tw}
\subsection{$\mathcal{K}$-structures and $\mathcal{S}$-structures}
In this section we recall some of the work from \cite{B} (see also \cite{Yano,Yano2}). We start with some definitions:
\begin{defi} An $f$-structure on a manifold $M^{2n+s}$ is a $(1,1)$ tensor field satisfying $f^3+f=0$.
\end{defi}
As mentioned in the introduction the existence of such a structure is equivalent to a reduction of the structural group of $M$ to $U(n)\times O(s)$. Throughout the paper we will use the notation "$M^{2n+s}$" to indicate that $M$ is a $(2n+s)$-dimensional manifold endowed with an $f$-structure with $rank(f)=2n$.
\begin{defi} We say that an $f$-structure on $M^{2n+s}$ has complemented frames if there are vector fields $\xi_i$ for  together with $1$-forms $\eta_i$ for $i\in\{1,...,s\}$ satisfying:
$$\eta_i(\xi_j)=\delta_{ij},\quad f\xi_i=0,\quad \eta_i\circ f=0,$$
$$f^2=-I+\sum_{k=1}^s\xi_k\otimes\eta_k$$
for all $i,j\in\{1,...,s\}$.
\end{defi}
\begin{rem} We list some of the immediate properties of $f$-structures with complemented frames:
\begin{enumerate}
\item $Ker(f)$ is equall to the bundle $<\xi_1,...,\xi_s>$ spammed pointwise by $\{\xi_1,...,\xi_s\}$. In particular this gives a partition $TM=Im(f)\oplus <\xi_1,...,\xi_s>$.
\item A complemented frame admits a compatible Riemannian metric $g$, i.e. such that:
$$g(X,Y)=g(fX,fY)+\sum_{k=1}^s \eta_k(X)\eta_k(Y).$$
with respect to such a metric the forms $\eta_i$ are dual to the corresponding vector fields $\xi_i$.
\item A compatible metric allows us to define a $2$-form $F(\bullet,\bullet):=g(\bullet,f\bullet)$ which is non-degenerate on $Im(f)$. It is easy to see that $\eta_1...\eta_s F^n\neq 0$. This implies that $M$ is orientable. Thorughout the rest of the paper we will consider $M$ with the orientation induced by the above $(2n+s)$-form
\end{enumerate}
\end{rem}
We now specialize to $\mathcal{K}$-manifolds introduced in \cite{B}. However, our main results can be applied to a slightly more general class of manifolds which is more natural to define beforehand:
\begin{defi} A manifold $M^{2n+s}$ together with an $f$-structure with a choosen complemented frame $\{\xi_1,...,\xi_s,\eta_1,...,\eta_s\}$ and a compatible Riemannian metric $g$ is an almost $\mathcal{K}$-manifold if:
\begin{enumerate}
\item The $2$-form $F(\bullet,\bullet):=g(\bullet,f\bullet)$ is closed.
\item The vector fields $\{\xi_1,...,\xi_s\}$ are Killing.
\end{enumerate}
If in addition the above set of data satisfies the equation:
$$ [f,f]+\sum_{k=1}^{s} \xi_k\otimes d\eta_k=0,$$
where $[f,f]$ is the Nijenhuis tensor of $f$, then the almost $\mathcal{K}$-structure is said to be integrable. An integrable almost $\mathcal{K}$-structure is also called a $\mathcal{K}$-structure (cf.\cite{B}).
\end{defi}

\begin{prop} Let $M^{2n+s}$ be an almost $\mathcal{K}$-manifold. Then $Ker(f)$ is involutive and hence induces a foliation.
\begin{proof} It follows from the definition that $Ker(f)$ is equall to the kernel of the closed 2-form $F(X,Y):=g(X,f Y)$ (treated as a map from $TM$ to $T^*M$). Hence, we get for any vector field $X$ the following equalities:
\begin{eqnarray*}
0=dF(\xi_i,\xi_j,X)&=&\mathcal{L}_{\xi_i}(F(\xi_j,X))-\mathcal{L}_{\xi_j}(F(\xi_i,X))+\mathcal{L}_{X}(F(\xi_i,\xi_j))
\\ &&-F([\xi_i,\xi_j], X)+F([\xi_i,X], \xi_j)-F([\xi_j,X], \xi_i)
\\ &=&-F([\xi_i,\xi_j], X)
\end{eqnarray*}
Which means that $[\xi_i,\xi_j]$ is again in the kernel of $F$ and hence the involutivity of $Ker(f)$ follows.
\end{proof}
\end{prop}
Moreover, it is clear that the foliation above is Riemannian (with $g(f\bullet,f\bullet)$) and transversely symplectic (with the $2$-form $F$), while the integrability condition implies that the foliation is also transversely holomorphic (and hence transversely K\"ahler). An almost $\mathcal{K}$-structure not only defines the above foliation with its additional structures but also connects the transverse geometry of the foliation to that of the entire manifold. One instance of this is the following proposition which we will use later on in some of our applications.
\begin{prop}\label{star} Let $M^{2n+s}$ be an almost $\mathcal{K}$-manifold and let $i=(i_1,...,i_k)$ be an ordered subset of $\{1,...,s\}$ with complement $j=(j_1,...,j_{s-k})$. Then the following relation between the hodge star operator $*$ and the basic hodge star operator $*_{b}$ holds for any transverse $r$-form $\alpha$ (i.e. $i_{\xi_l}\alpha=0$):
$$*(\eta_{i_1}...\eta_{i_k}\alpha)=(-1)^{sign(i_1,...,i_k,j_1,...,j_{s-k})+(s-k)r}\eta_{j_1}...\eta_{j_{s-k}}*_b\alpha$$
\end{prop}

Among $\mathcal{K}$-structure the notions of $\mathcal{S}$-structure and $\mathcal{C}$-structures was given special attention in \cite{B} due to them being the proper generalization of Sasakian and quasi-Sasakian with $d\eta=0$ cases to the above setting and as such exhibit analogous curvature properties. Again we introduce these structures proceeded by their "almost structure" counterpart:
\begin{defi} An almost $\mathcal{K}$-structure on a manifold $M^{2n+s}$:
\begin{enumerate}
\item is called an almost $\mathcal{S}$-structure if $d\eta_i=F$ for all $i\in\{1,...,s\}$.
\item is called an almost $\mathcal{C}$-structure if $d\eta_i=0$ for all $i\in\{1,...,s\}$.
\end{enumerate}
Moreover, if the underlying almost $\mathcal{K}$-structure of an almost $\mathcal{S}$-structure (resp. almost $\mathcal{C}$-structure) is integrable, then it is called a $\mathcal{S}$-structure (resp. $\mathcal{C}$-structure).
\end{defi}
We finish this section by reiterating for the readers convienience the analogy between the above structures and their low dimensional counterparts.
$$\begin{tabular}{|c|c|c|}\hline s=0 & s=1 & general\\\hline - & Quasi-K-contact & Almost $\mathcal{K}$\\- & Quasi-Sasaki & $\mathcal{K}$\\
Almost K\"ahler &K-contact& Almost $\mathcal{S}$\\
K\"ahler & Sasaki &  $\mathcal{S}$\\
 - & Quasi-K-contact with $d\eta=0$& Almost $\mathcal{C}$\\- & Quasi-Sasaki with $d\eta=0$ & $\mathcal{C}$
\\\hline
\end{tabular}$$
\begin{rem} It is also worth noting that our almost $\mathcal{S}$-manifolds are already present in the literature as $f$-$K$-contact manifolds. However, we choose to stick to our terminology as it seems more appropriate when considering such manifolds along with non-integrable versions of $\mathcal{K}$-manifolds and $\mathcal{C}$-manifolds.
\end{rem}
\section{The spectral sequence of invariant forms on almost $\mathcal{K}$-manifolds}\label{Tor}
Here we describe a canonical torus action on certain almost $\mathcal{K}$-manifolds and use it to define a spectral sequence used in further chapters. We start with the following proposition:
\begin{prop}\label{com} Let $M^{2n+s}$ be an almost $\mathcal{K}$-manifold. such that for each $l\in\{1,...,s\}$ the form $d\eta_l$ is basic. Then for each $i,j\in\{1,...,s\}$ the equality $[\xi_i,\xi_j]=0$ holds.
\begin{proof}
This follows from the computation:
\begin{eqnarray*}
0=d\eta_l(\xi_i,\xi_j)=\mathcal{L}_{\xi_i}\eta_l(\xi_j)-\mathcal{L}_{\xi_j}\eta_l(\xi_i)-\eta_{l}([\xi_i,\xi_j])=-\eta_l([\xi_i,\xi_j]).
\end{eqnarray*}
Now since $Ker(f)$ is involutive it follows that if $[\xi_i,\xi_j]\neq 0$ then there exists some $l$ such that $\eta_l([\xi_i,\xi_j])\neq 0$ which provides the desired contradiction with the computation above.
\end{proof}
\end{prop}
\begin{rem} Let us briefly note that in particular almost $\mathcal{S}$-manifolds and almost $\mathcal{C}$-manifolds satisfy the asumptions of the above proposition. Hence, this proposition as well as the remainder of this section can be applied to them.
\end{rem}
This has the following important corollary:
\begin{cor} Let $M^{2n+s}$ be a compact almost $\mathcal{K}$-manifold. such that for each $i\in\{1,...,s\}$ the form $d\eta_i$ is basic. Let $G\subset Diff(M)$ be the group whose Lie algebra is $<\xi_1,...,\xi_s>\subset \Gamma(TM)$. Then the closure of $G$ is a Torus in the group $Isom(M)$ of isometries on $M$.
\begin{proof} Let us first note that since the vector fields $\xi_i$ are killing we have the inclusion $G\subset Isom(M)$ which is known to be a finitely dimensional compact Lie group. Moreover, since by Proposition \ref{com} $G$ is abelian its closure is a compact abelian group and hence a torus.
\end{proof}
\end{cor}
The next step is to classify forms on $M$ which are invariant under the action of $\overline{G}$.
\begin{prop}Let $M^{2n+s}$ be a compact almost $\mathcal{K}$-manifold such that for each $i\in\{1,...,s\}$ the form $d\eta_i$ is basic. Then the following conditions are equivalent:
\begin{enumerate}
\item $\alpha$ is a $\overline{G}$-invariant form on $M$.
\item $\alpha=\alpha_0+\sum\limits_{k=1}^s\sum\limits_{1\leq i_1<...<i_k\leq s}\eta_{i_1}...\eta_{i_k}\alpha_{i_1,...,i_k}$, where $\alpha_0$ and $\alpha_{i_1,...,i_k}$ are basic for all indices $1\leq i_1<...<i_k\leq s$.
\end{enumerate}
\begin{proof} Assume that the second condition is true. Then it can be easilly computed that for any $\xi_j$ the equality $\mathcal{L}_{\xi_j}\alpha=0$ holds. Which in turn implies that $\alpha$ is $\overline{G}$-invariant.
\newline\indent Now let us write the invariant form $\alpha$ as:
$$\alpha=\alpha_0+\sum\limits_{k=1}^s\sum\limits_{1\leq i_1<...<i_k\leq s}\eta_{i_1}...\eta_{i_k}\alpha_{i_1,...,i_k},$$
where $\alpha_{i_1,...,i_k}$ are transverse for all indices $1\leq i_1<...<i_k\leq s$. Due to the well known formula:
$$\mathcal{L}_{X}i_{Y}-i_{Y}\mathcal{L}_{X}=i_{[X,Y]},$$
we get that $i_{\xi_i}$ and $\mathcal{L}_{\xi_j}$ commute for $i,j\in\{1,...,s\}$ (using Proposition \ref{com}). We shall now prove that the forms $\alpha_0$ and $\alpha_{i_1,...,i_k}$ are basic by reverse induction on the number of indices. Hence, we start by proving that $\alpha_{1,...,s}$ is basic. Since $\alpha$ is harmonic and the vector fields $\xi_i$ are Killing we have for any $i\in\{1,...,s\}$ the following equalities:
$$0=\mathcal{L}_{\xi_i}\alpha=i_{\xi_s}i_{\xi_{s-1}}...i_{\xi_1}\mathcal{L}_{\xi_i}\alpha=\mathcal{L}_{\xi_i}i_{\xi_{s}}i_{\xi_{s-1}}...i_{\xi_1}\alpha=\mathcal{L}_{\xi_i}\alpha_{1,...,s}.$$
Which proves that $\alpha_{1,...,s}$ is basic.
\newline\indent For the induction step let us assume that all the $\alpha_{i_1,...,i_{k}}$ for $s\geq k>K$ are basic. We shall show that all $\alpha_{i_1,...,i_{K}}$ are basic as well. Using the assumption we get for any $i\in\{1,...,s\}$ the following equalities:
$$0=\mathcal{L}_{\xi_i}\alpha=i_{\xi_{i_K}}i_{\xi_{i_{K-1}}}...i_{\xi_{i_1}}\mathcal{L}_{\xi_i}\alpha=\mathcal{L}_{\xi_i}i_{\xi_{i_K}}i_{\xi_{i_{K-1}}}...i_{\xi_{i_1}}\alpha=\mathcal{L}_{\xi_i}\alpha_{i_1,...,i_K}.$$
Which proves that $\alpha_{i_1,...,i_K}$ are basic for any set of indices $1\leq i_1<...<i_k\leq s$.
\end{proof}
\end{prop}
\begin{rem} Note that the induction assumption is used to pass to the final equality as it implies that all the terms with a greater number of indices then $K$ vanish under $\mathcal{L}_{\xi_i}$ as:
$$\mathcal{L}_{\xi_i}\eta_{j_1}...\eta_{j_k}\alpha_{j_1,...,j_k}=\eta_{j_1}...\eta_{j_k}\mathcal{L}_{\xi_i}\alpha_{j_1,...,j_k}=0.$$
The first equality is due to the fact that $\mathcal{L}_{\xi_i}\eta_j=i_{\xi_i}d\eta_j+d(i_{\xi_i}\eta_j)=0$.
\end{rem}

Finally, we note that similarilly as for the spectral sequence of a Riemannian foliation we have a filtration of the cochain complex of invariant forms $\Omega^r_{\overline{G}}(M)$ given by:
$$F^k_{\mathcal{F}}\Omega^r_{\overline{G}}(M):=\{\alpha\in\Omega^r_{\overline{G}}(M)\text{ }|\text{ } i_{X_{r-k+1}}...i_{X_1}\alpha=0,\text{ for } X_1,...,X_{r-k+1}\in\Gamma(T\mathcal{F})\}.$$
Hence, via known theory from homological algebra we can cosntruct a spectral sequence as follows:
\begin{enumerate}
\item The $0$-th page is given by $E_0^{p,q}=F_{\mathcal{F}}^{p}\Omega^{p+q}_{\overline{G}}(M)\slash F_{\mathcal{F}}^{p+1}\Omega^{p+q}_{\overline{G}}(M)$ and $d^{p,q}_0:E_0^{p,q}\to E^{p,q+1}_0$ is simply the morphism induced by $d$.
\item The $r$-th page is given inductively by:
$$E_{r}^{p,q}:=Ker(d^{p,q}_{r-1})\slash Im(d^{p,q}_{r-1})=\frac{\{\alpha\in F^{p}_{\mathcal{F}}\Omega^{p+q}_{\overline{G}}(M)\text{ }|\text{ } d\alpha\in F_{\mathcal{F}}^{p+r}\Omega^{p+q+1}_{overline{G}}(M)\}}{F^{p+1}_{\mathcal{F}}\Omega^{p+q}_{\overline{G}}(M)+d(F_{\mathcal{F}}^{p-r+1}\Omega_{\overline{G}}^{p+q-1}(M))}$$
\item The $r$-th coboundary operator $d_r:E_r^{p,q}\to E^{p+r,q-r+1}_r$ is again just the map induced by $d$ (due to the description of the $r$-th page this has the target specified above and is well defined).
\end{enumerate}
Furthermore, since the filtration is bounded this spectral sequence converges and its final page is isomorphic to the cohomology of the cochain complex $\Omega^{r}_{\overline{G}}(M)$ known to be isomorphic to the de Rham cohomology of $M$. We call this spectral sequence the spectral sequence of invariant forms and denote it by $E_{r}^{p,q}$ throughout the rest of the paper.
\begin{tw}\label{2nd} Let $M^{2n+s}$ be a compact almost $\mathcal{S}$-manifold such that for each $i\in\{1,...,s\}$ the form $d\eta_i$ is basic. Then:
$$E_2^{p,q}\cong \bigwedge\text{}_{H^{p}(M\slash\mathcal{F})}^q<\eta_1,...,\eta_s>:=H^{p}(M\slash\mathcal{F})\otimes\bigwedge\text{} ^q<\eta_1,...,\eta_s>.$$
\begin{proof} Since the operator $d$ takes basic forms to basic forms and $d\eta_i$ is basic for all $i\in\{1,...,s\}$ it is easy to see that $d_0$ is in fact equall to the zero operator. Hence, the first page is isomorphic to the $0$-th page.
\newline\indent On the first page by the same observation the operator $d_1$ is just the application of $d$ to the transverse part of the form (since applying $d$ to $\bigwedge ^q<\eta_1,...,\eta_s>$ decrease $q$). Hence, the second page is just $H^{p}(M\slash\mathcal{F})\otimes\bigwedge ^q<\eta_1,...,\eta_s>$.
\end{proof}
\end{tw}
\begin{rem} \begin{enumerate} 
\item We note that the merit of considering invariant forms is already visible in this computation since a similar result for the spectral sequence of the Riemannian foliation is not known. In fact, it is far from trivial to even proof that the second page of this spectral sequence is finitely dimensional (cf. \cite{ss1}). On a more down to earth level the major simplification comes from the triviality of $d_0$ in the spectral sequence of invariant forms.
\item It is interesting to note that for $\mathcal{S}$-manifolds the above description of $E^{p,q}_2$ coincides (disregarding the coboundary operators) with the "almost formal" models from \cite{Vais}.
\item For the sake of brieviety the notation:
$$\bigwedge\text{}_{V}^q<\eta_1,...,\eta_s>:=V\otimes\bigwedge\text{} ^q<\eta_1,...,\eta_s>,$$
introduced in the above theorem, shall be used throughout the article. We shall also use its following variations:
$$\bigwedge\text{}_{V}^{\bullet}<\eta_1,...,\eta_s>:=V\otimes\bigwedge\text{} ^{\bullet}<\eta_1,...,\eta_s>,$$
$$\overline{\bigwedge}\text{}_{V}^{\bullet}<\eta_1,...,\eta_s>:=\{\alpha\in V\otimes\bigwedge\text{} ^{\bullet}<\eta_1,...,\eta_s>\text{ }|\text{ } \pi_{V\otimes 1}\alpha=0\},$$
where $\pi_{V\otimes 1}$ is the obvious projection onto $V\otimes 1\subset \bigwedge\text{}_{V}^{\bullet}<\eta_1,...,\eta_s>.$
\end{enumerate}
\end{rem}
We also wish to mention the following consequence of the above discussion which will be used throughout the paper in order to omit the homological orientability assumption throughout the rest of the article:
\begin{prop} Let $M^{2n+s}$ be a compact almost $\mathcal{S}$-manifold such that for each $i\in\{1,...,s\}$ the form $d\eta_i$ is basic. Then the foliation induced by $Ker(f)$ is homologically orientable.
\begin{proof} It is well known (cf. \cite{E1}) that the top basic cohomology of a Riemannian foliation on a compact manifold is either $0$ or $\mathbb{R}$. In this case it cannot be $0$ since then we could compute from the above spectral sequence that $H^{2n+s}_{dR}(M)\cong E^{2n,s}_2=0$ which is a contradiction with the orientability of $M$. Hence, the top basic cohomology is isomorphic to $\mathbb{R}$ which means that the foliation is homologically orientable.
\end{proof}
\end{prop}
\section{Main results for almost $\mathcal{S}$-manifolds}\label{MainS}
In this section we prove our main results for compact almost $\mathcal{S}$-manifolds. We start by computing $E^{p,q}_3$ for almost $\mathcal{S}$-manifolds satisfying the transverse hard Lefschetz property. Firstly, we compute the kernel of $d_2$.
\begin{lem}\label{ker} Under the above assumptions we have:
$$Ker(d_2)=(\bigwedge\text{}_{H^{\bullet}(M\slash\mathcal{F})}^{\bullet} <\eta_1-\eta_{2},...,\eta_1-\eta_{s}> +\overline{\bigwedge}\text{}_{Ker(L)}^{\bullet} <\eta_{1},...,\eta_{s}>).$$
\begin{proof} Firstly, let us note that $Ker(d^{p,0}_2)$ is simply $H^{p}(M\slash\mathcal{F})$ . For $q>0$, by Theorem \ref{2nd} we can write an element $[\alpha]$ of $E_2^{p,q}$ as:
$$[\alpha]=\sum\limits_{1\leq i_1<...<i_q\leq s}\eta_{i_1}...\eta_{i_q}[\alpha_{i_1,...,i_q}],$$
for some basic forms $[\alpha_{i_1,...,i_q}]$.
Applying $d_2$ to this element gives:
$$d_2[\alpha]=\sum\limits_{1\leq i_1<...<i_{q-1}\leq s}\eta_{i_1}...\eta_{i_{q-1}}L(\sum\limits_{j=1}^s[\alpha_{i_1,...,i_{q-1},j}])$$
Where we understand $\alpha_{i_1,...,i_{q-1},j}$ to be equal to zero if $j\in\{i_1,...,i_{q-1}\}$ and to be equal to $sign(i_1,...,i_{q-1},j)\alpha_{i_1,...,j,...,i_{q-1}}$ otherwise (here $j$ is on the correct position so that the indices are in an increasing order). This implies that $\overline{\bigwedge}_{Ker(L)}^{\bullet} <\eta_{1},...,\eta_{s}>$ is in fact contained in $Ker(d^{p,q}_2)$.
\newline\indent Here we split the consideration into two cases $p< n$ and $p\geq n$.
\newline\indent For the first case, $L$ is a monomorphism and the elements $\sum\limits_{j=1}^s[\alpha_{i_1,...,i_{q-1},j}]$ have to be trivial for $\alpha$ to be an element of $Ker(d_2^{p,q})$. However, by assigning $\eta_{i_1}...\eta_{i_k}$ to the simplex $[i_1,...,i_k]$ we get a commutative diagram:
$$\begin{tikzcd}
 \tilde{C}_{q}(\Delta^s;H^{p}(M\slash\mathcal{F}))\arrow[r]\arrow[d,"\partial^q"]& \bigwedge\text{}_{H^{p}(M\slash\mathcal{F})}^{q} <\eta_1,...,\eta_{s}>\arrow[d,"d_2^{p,q}"]
\\ \tilde{C}_{q-1}(\Delta^s;H^{p}(M\slash\mathcal{F}))\arrow[r]& \bigwedge\text{}_{LH^{p}(M\slash\mathcal{F})}^{q-1} <\eta_1,...,\eta_{s}>
\end{tikzcd}.$$
The horizontal arrows in this diagram are isomorphisms and hence they induce an isomorphism on the kernels of the vertical arrows. But $Ker(\partial^q)=Im(\partial^{q+1})$. While for $Im(\partial^{q+1})$ we can easilly determine the generators as the images of the simplices generating $\tilde{C}^{q+1}$. These in turn correspond to the elements of the form $(\eta_{i_1}-\eta_{i_2})...(\eta_{i_1}-\eta_{i_q})\alpha$ for $\alpha\in H^{p}(M\mathcal{F})$.
\newline\indent For the second case, Theorem \ref{prim1} implies that, $Ker(L)$ is equall to $L^{p-n}\mathcal{P}H^{2n-p}(M\slash\mathcal{F})$.
Hence, in what follows it suffices to consider classes from $L^{p-n+1}H^{2n-p-2}(M\slash\mathcal{F})$ on which $L$ is monomorphic. From here the same argument as in the first case can be conducted with the coefficients changed to $L^{p-n+1}H^{2n-p-2}(M\slash\mathcal{F})$.
\end{proof}
\end{lem}
With this we are ready to prove our main result concerning almost $\mathcal{K}$-manifolds:
\begin{tw}\label{mainS} Let $M^{2n+s}$ be a compact almost $\mathcal{K}$-manifold satisfying the transverse hard Lefschetz property. Then:
$$H^{\bullet}_{dR}(M) \cong \bigwedge\text{}^{\bullet}_{\mathcal{P}H^{\bullet}(M\slash\mathcal{F})}<\eta_1-\eta_2,...,\eta_1-\eta_s>\oplus\eta_1\bigwedge\text{}^{\bullet}_{Ker(L)}<\eta_2,...,\eta_s>.$$
\begin{proof} Firstly, let us note that the image of $d^{p,1}_2$ is equall to the image of $L$ and hence $E^{p,0}_3\cong\mathcal{P}H^{p}(M\slash\mathcal{F})$.
\newline\indent Secondly, for $p<n$ we again have that $L$ is a monomorphism and hence we have the identification:
$$\begin{tikzcd}
 \tilde{C}_{q+1}(\Delta^s;H^{p-2}(M\slash\mathcal{F}))\arrow[r]\arrow[d,"\partial^{q+1}"]& \bigwedge\text{}_{H^{p-2}(M\slash\mathcal{F})}^{q+1} <\eta_1,...,\eta_{s}>\arrow[d,"d_2^{p-2,q+1}"]
\\ \tilde{C}_{q}(\Delta^s;H^{p-2}(M\slash\mathcal{F}))\arrow[r]& \bigwedge\text{}_{LH^{p-2}(M\slash\mathcal{F})}^{q} <\eta_1,...,\eta_{s}>
\end{tikzcd}.$$
This implies (similarilly as in the previous proof) that the image consists of elements of the form $(\eta_{i_1}-\eta_{i_2})...(\eta_{i_1}-\eta_{i_q})L\alpha$ for $\alpha\in H^{p-2}(M\mathcal{F})$. With this (together with Lemma \ref{ker}) we conclude that in this range of $p$ we have:
$$E_3^{p,q}\cong\bigwedge\text{}^{\bullet}_{\mathcal{P}H^{\bullet}(M\slash\mathcal{F})}<\eta_1-\eta_2,...,\eta_1-\eta_s>$$.
\newline\indent Thirdly, we treat the case $p\geq n$. Here it is crucial to note that due to its form $d^{p,q}_2$ preserves the basic Lefschetz decomposition of each element of $E^{p,q}_2$. In particular, this allows us to consider seperately $d^{p-2,q+1}_2$ on $L^{p-n+2}H^{2n-p-4}(M\slash\mathcal{F})$ and $L^{p-n+1}\mathcal{P}H^{2n-p-2}(M\slash\mathcal{F})$.
\newline\indent For $L^{p-n+2}H^{2n-p-4}(M\slash\mathcal{F})$, a similar diagram as in the previous case (but with coefficients in $L^{p-n+2}H^{2n-p-4}(M\slash\mathcal{F})$) will allow us to compute the image of $d^{p-2,q+1}_2$ to consist of elements of the form $(\eta_{i_1}-\eta_{i_2})...(\eta_{i_1}-\eta_{i_q})L\alpha$ for $\alpha\in L^{p-n+2}H^{2n-p-4}(M\slash\mathcal{F})$. But since in the given range of $p$ the morphism $L$ is an epimorphism this indicates that for such $p$ the algebra $\bigwedge\text{}_{H^{\bullet}(M\slash\mathcal{F})\slash Ker(L)}^{\bullet} <\eta_1-\eta_{2},...,\eta_1-\eta_{s}>$ trivializes on $E^{p,q}_3$.
\newline\indent For $L^{p-n+1}\mathcal{P}H^{2n-p-2}(M\slash\mathcal{F})$, we again use the same method to compute that the image of $d^{p-2,q+1}_2$ consist of elements of the form $(\eta_{i_1}-\eta_{i_2})...(\eta_{i_1}-\eta_{i_q})L\alpha$ for $\alpha\in L^{p-n+1}\mathcal{P}H^{2n-p-2}(M\slash\mathcal{F})$. This by the basic Lefschetz decomposition can be written alternatively as $(\eta_{i_1}-\eta_{i_2})...(\eta_{i_1}-\eta_{i_q})\alpha$ for $\alpha\in Ker(L)$. Combining this with the previous paragraph we get that in this range of $p$ the following holds:
$$E^{p,q}_3\cong \overline{\bigwedge}\text{}_{Ker(L)}^{\bullet} <\eta_{1},...,\eta_{s}>\slash\overline{\bigwedge}\text{}_{Ker(L)}^{\bullet} <\eta_1-\eta_{2},...,\eta_1-\eta_{s}>.$$
This in turn can be easilly computed to be isomorphic to the complement of $\overline{\bigwedge}\text{}_{Ker(L)}^{\bullet} <\eta_1-\eta_{2},...,\eta_1-\eta_{s}>$ given for example by $\eta_1\overline{\bigwedge}\text{}^{\bullet}_{Ker(L)}<\eta_2,...,\eta_s>.$
\newline\indent Finally, we note that by the basic Lefschetz decomposition we can take the basic components of the representatives of the classes from $E_3^{p,q}$ to be either primitive closed forms or closed forms from $Ker(L)$ which by the construction of the spectral sequence proves that the spectral sequence degenerates at the 3rd page. Moreover, by pinpointing the representatives in such a way we can conclude that they also represent the cohomology classes in $H^{\bullet}_{dR}(M)$.
\end{proof}
\end{tw}
\begin{rem}\label{GL} We note that a different proof of Theorem \ref{mainS} can be conducted by combining some recent results from \cite{GL} with methods from \cite{lin}. More precisely, Proposition 4.4 from \cite{GL} gives a Gysin-like long exact sequence connecting the basic cohomology of the foliation $\mathcal{F}_{s-1}$ spaned by $\{\xi_2,...,\xi_{s}\}$ and the basic cohomology of $\mathcal{F}$. By analyzing it (with the use of the basic hard Lefschetz property) similarilly as in \cite{lin} we get:
$$H^{\bullet}(M\slash\mathcal{F}_{s-1})\cong \mathcal{P}H^{\bullet}\oplus \eta_1 Ker(L).$$
By using now Theorem 4.5 from \cite{GL}, which relates $H^{\bullet}(M\slash\mathcal{F}_{s-1})$ to $H^{\bullet}_{dR}(M)$ one arrives at the conclusion of Theorem \ref{mainS}.
\end{rem}
An immediate consequence of Theorem \ref{mainS} is the following important result on topological invariance of basic cohomology.  
\begin{cor}\label{TopInvS} Let $M^{2n+s}_1$ and $M^{2n+s}_2$ be compact almost $\mathcal{S}$-manifolds satisfying the hard Lefschetz property which are homeomorphic. Then the basic cohomologies of the corresponding foliations are isomorphic.
\begin{proof} The above description can be used to compute the primitive basic cohomology of both the foliations $\mathcal{F}_l$ on $M_l$ for $l\in\{1,2\}$. More precisely it can be done inductively by the formula:

$$dim(\mathcal{P}H^{k}(M_l\slash\mathcal{F}_l))=dim(H^{k}_{dR}(M)_l)-\sum\limits^{k-1}_{i=0} {s-1\choose k-i} dim(\mathcal{P}H^i(M_l\slash\mathcal{F}_l))$$
with the convention that ${s-1\choose k-i}=0$ if $k-i>s-1$. This implies that:
$$dim(\mathcal{P}H^{k}(M_1\slash\mathcal{F}_1))=dim(\mathcal{P}H^{k}(M_2\slash\mathcal{F}_2)),$$
for all $0\leq k\leq n$. But these dimensions are enough to compute the dimensions of the basic cohomology of $\mathcal{F}$ using the Lefschetz decomposition. Hence:
$$dim(H^{k}(M_1\slash\mathcal{F}_1))=dim(H^{k}(M_2\slash\mathcal{F}_2)),$$
for all $k\in\mathbb{N}$ which in turn implies that the basic cohomologies of $\mathcal{F}_1$ and $\mathcal{F}_2$ are isomorphic.
\end{proof}  
\end{cor}
\begin{rem} Unfortunately, similarilly as in the Sasakian (and K-contact) case this method does not produce any cannonical isomorphism between the basic cohomologies of $\mathcal{F}_1$ and $\mathcal{F}_2$.
\end{rem}

\section{Main result for almost $\mathcal{C}$-manifolds}\label{MainC}
\begin{tw}\label{mainC} Let $M^{2n+s}$ be an almost $\mathcal{C}$-manifold. Then:
$$H^{k}_{dR}(M)\cong \bigoplus\limits_{p+q=k} H^{p}(M\slash\mathcal{F})\otimes\bigwedge\text{} ^q<\eta_1,...,\eta_s>.$$
\begin{proof} In this case the operator $d$ itself goes from $\Omega^{p}(M\slash\mathcal{F})\otimes\bigwedge ^q<\eta_1,...,\eta_s>$ to $\Omega^{p+1}(M\slash\mathcal{F})\otimes\bigwedge ^q<\eta_1,...,\eta_s>$. Hence, the spectral sequence degenerates at the second page which together with Theorem \ref{2nd} implies the thesis. 
\end{proof}
\end{tw}
It is worth noting that in this case the transverse hard Lefschetz property is not needed. In particular, a similar result holds for quasi-Sasakian manifolds with $d\eta=0$.Similatily as in the $\mathcal{S}$-manifold case this also implies that the basic cohomology of almost $\mathcal{C}$-manifolds are a topological invariant.

\begin{cor}\label{TopInvC}
  Let $M^{2n+s}_1$ and $M^{2n+s}_2$ be compact almost $\mathcal{C}$-manifolds which are homeomorphic. Then the basic cohomologies of the induced foliations are isomorphic. 
\begin{proof} The above description can be used to compute the basic cohomology of both the foliations $\mathcal{F}_l$ on $M_l$ for $l\in\{1,2\}$. More precisely it can be done inductively by the formula:

$$dim(H^{k}(M_l\slash\mathcal{F}_l))=dim(H^{k}_{dR}(M)_l)-\sum\limits^{k-1}_{i=0} {s\choose k-i} dim(H^i(M_l\slash\mathcal{F}_l))$$
with the convention that ${s\choose k-i}=0$ if $k-i>s$. This implies that:
$$dim(H^{k}(M_1\slash\mathcal{F}_1))=dim(H^{k}(M_2\slash\mathcal{F}_2)),$$
for all $k\in\mathbb{N}$ which in turn implies that the basic cohomologies of $\mathcal{F}_1$ and $\mathcal{F}_2$ are isomorphic.
\end{proof}
\end{cor}
\begin{rem} As with the anologous result from the previous chapter this method does not produce any cannonical isomorphism between the basic cohomologies of $\mathcal{F}_1$ and $\mathcal{F}_2$.
\end{rem}
\begin{rem} While it seems doubtful that similar general results can be achieved for arbitrary almost $\mathcal{K}$-manifolds we feel that the above spectral sequence remains a good tool to find similar dependencies in cohomology on a case by case basis.
\end{rem}
\section{Application: Classification of Harmonic forms}\label{AppHF}
Here we provide a description of Harmonic forms on almost $\mathcal{C}$-manifolds and $\mathcal{S}$-manifolds based on their basic harmonic forms. Firstly, let us note the following:
\begin{rem}\label{remS}\begin{enumerate}
\item Theorem \ref{mainS} allows us to treat $\bigwedge\text{}^{\bullet}_{\mathcal{P}H^{\bullet}(M\slash\mathcal{F})}<[\eta_1-\eta_2],...,[\eta_1-\eta_s]>$ and $\eta_1\bigwedge\text{}^{\bullet}_{Ker(L)}<\eta_2,...,\eta_s>$ as submodules of $H^{\bullet}_{dR}(M)$ with representatives given respectively by linear combinations of elements of the form $(\eta_{1}-\eta_{i_1})...(\eta_1-\eta_{i_q})\alpha$ (where $\alpha$ is a closed basic primitive form) and $\eta_1\eta_{i_1}...\eta_{i_s}\alpha$ (where $\alpha\in Ker(L)$ is a closed basic form).
\item Theorem \ref{mainC} allows us to treat $[\eta_{i_1}]...[\eta_{i_q}]H^{p}(M\slash\mathcal{F})$ as a submodule of $H^{p+q}_{dR}(M)$ with representatives given by linear combinations of elements of the form $\eta_{i_1}...\eta_{i_q}\alpha$ (where $\alpha$ is a closed basic form).
\end{enumerate}
\end{rem}
Making the identification from the previous remark the following statements can be now made:
\begin{tw} Let $M^{2n+s}$ be a compact almost $\mathcal{C}$-manifold. A form $\alpha$ on $M$ is harmonic if and only if it is a linear combination of elements of the form $\eta_{i_1}...\eta_{i_q}\tilde{\alpha}$ such that $\tilde{\alpha}$ is basic harmonic.

\begin{proof} We start by noting that any element $\alpha$ of the given form is in fact harmonic (which immediately implies that the linear combination of such elements is also harmonic), since through straightforward computation (with the use of proposition \ref{star}) we can get $d\alpha=\delta\alpha=0$.
\newline\indent On the other hand given any harmonic form $\alpha$ by Theorem \ref{mainC} the cohomology class it represents splits into a sum of elements of the form:
$$[\eta_{i_1}]...[\eta_{i_q}][\tilde{\alpha}]\in [\eta_{i_1}]...[\eta_{i_q}]H^{p}(M\slash\mathcal{F})\subset H^{p+q}_{dR}(M).$$
Hence, we can write this harmonic form as the sum of the harmonic forms corresponding to such classes. The proof is now finished by noting that the forms $\eta_{i_1}...\eta_{i_q}\tilde{\alpha}$ where $\tilde{\alpha}$ is basic harmonic are the representatives of such classes (as this implies that they are basic harmonic by the previous paragraph).
\end{proof}
\end{tw}
\begin{tw}Let $M^{2n+s}$ be a compact $\mathcal{S}$-manifold. A form $\alpha$ on $M$ is harmonic if and only if it is a linear combination of elements of the form $(\eta_1-\eta_{i_1})...(\eta_1-\eta_{i_q})\tilde{\alpha}$ such that $\tilde{\alpha}$ is basic primitive harmonic and their duals (via the star operator).
\begin{proof} We start by noting that any element $\alpha$ of the given form is in fact harmonic (which immediately implies that the linear combination of such elements is also harmonic), since through straightforward computation (with the use of proposition \ref{star}) we can get $d\alpha=\delta\alpha=0$.
\newline\indent On the other hand given any harmonic form $\alpha$ by Theorem \ref{mainS} the cohomology class it represents splits into a sum of elements of one of the following forms:
$$[\eta_1-\eta_{i_1}]...[\eta_1-\eta_{i_q}][\tilde{\alpha}]\in [\eta_1-\eta_{i_1}]...[\eta_1-\eta_{i_q}]\mathcal{P}H^{p}(M\slash\mathcal{F})\subset H^{p+q}_{dR}(M),$$
$$[\eta_1\eta_{i_1}...\eta_{i_{q-1}}\tilde{\alpha}]\in \eta_1\eta_{i_1}...\eta_{i_{q-1}}Ker^p(L)\subset H^{p+q}_{dR}(M).$$
Hence, it suffices to find the harmonic representatives of such classes.
\newline\indent For the classes in $[\eta_1-\eta_{i_1}]...[\eta_1-\eta_{i_q}]\mathcal{P}H^{p}(M\slash\mathcal{F})$ let us first note that by the Lefschetz decomposition of harmonic forms the basic harmonic representative of a class $[\tilde{\alpha}]\in\mathcal{P}H^p(M\slash\mathcal{F})$ is indeed primitive. Now it suffices to note that the class $[\eta_1-\eta_{i_1}]...[\eta_1-\eta_{i_q}][\tilde{\alpha}]$ is indeed represented by $(\eta_1-\eta_{i_1})...(\eta_1-\eta_{i_q})\tilde{\alpha}$ such that $\tilde{\alpha}$ is the basic harmonic representative of the basic class $[\tilde{\alpha}]$.
\newline\indent For the classes in $\eta_1\bigwedge\text{}^{\bullet}_{Ker(L)}<\eta_2,...,\eta_s>$ we shall prove that each such class is represented by a linear combination of duals of harmonic forms representing a class in $[\eta_1-\eta_{i_1}]...[\eta_1-\eta_{i_q}]\mathcal{P}H^{p}(M\slash\mathcal{F})$. Let us start by taking the harmonic representative $(\eta_1-\eta_{i_1})...(\eta_1-\eta_{i_q})\tilde{\alpha}$ of a class:
$$[\eta_1-\eta_{i_1}]...[\eta_1-\eta_{i_q}][\tilde{\alpha}]\in [\eta_1-\eta_{i_1}]...[\eta_1-\eta_{i_q}]\mathcal{P}H^{p}(M\slash\mathcal{F})\subset H^{p+q}_{dR}(M).$$
We can see by proposition \ref{star} (and the fact that the basic star operator takes basic primitive harmonic forms to harmonic forms in $Ker(L)$) that its dual is an element of $\overline{\bigwedge}\text{}_{Ker(L)}^{\bullet} <\eta_{1},...,\eta_{s}>$. Hence, by Lemma \ref{ker} (and the proof of Theorem \ref{mainS}) it represents a class in:
$$\overline{\bigwedge}\text{}_{Ker(L)}^{\bullet} <\eta_{1},...,\eta_{s}>\slash\overline{\bigwedge}\text{}_{Ker(L)}^{\bullet} <\eta_1-\eta_{2},...,\eta_1-\eta_{s}>\cong (\eta_1\bigwedge\text{}^{\bullet}_{Ker(L)}<\eta_2,...,\eta_s>)\subset H^{p+q}_{dR}(M).$$  
This shows that $*$ induces a morphism:
$$[*]: \bigwedge\text{}^{\bullet}_{\mathcal{P}H^{\bullet}(M\slash\mathcal{F})}<[\eta_1-\eta_2],...,[\eta_1-\eta_s]>\to \eta_1\bigwedge\text{}^{\bullet}_{Ker(L)}<\eta_2,...,\eta_s>,$$
by acting on the harmonic representatives. Moreover, due to well known Hodge theoretic results the star operator induces a bijection on harmonic forms and hence the morphism $[*]$ is a monomorphism (since it is a restriction of the star operator composed with the morphism induced by the inclusion of harmonic forms). From this we deduce that this morphism is in fact an isomorphism, since it is a monomorphism between vector spaces of the same (finite) dimension. This together with the fact that the star operator preserves harmonic forms shows that the harmonic representatives of the classes from $\eta_1\bigwedge\text{}^{\bullet}_{Ker(L)}<\eta_2,...,\eta_s>$ are precisely the duals of the harmonic forms given as linear combinations of elements of the form $(\eta_1-\eta_{i_1})...(\eta_1-\eta_{i_q})\tilde{\alpha}$ such that $\tilde{\alpha}$ is basic primitive harmonic.
\end{proof}
\end{tw}
\begin{rem} Again similarilly as in the previous section we believe that the methods used above can be useful in computing harmonic forms of general $\mathcal{K}$-manifolds on a case by case basis.
\end{rem}

\section{Application: Stability of basic Betti and Hodge numbers}\label{AppStab}
Here we wish to study the behaviour of the basic cohomology of an almost $\mathcal{K}$-manifold under deformations. Let us start by making the notion of a deformation of an almost $\mathcal{K}$-manifold precise.
\begin{defi}\label{defdef} Let $M^{2n+s}$ be a compact almost $\mathcal{K}$-manifold. A deformation $\{M_t\}_{t\in [0,1]}$ of $M$ consists of the following data:
\begin{enumerate}
\item a $(0,2)$-tensor $g$ on $M\times [0,1]$ such that its restriction $g_t$ to each $M_t:=M\times \{t\}$ is a Riemannian metric and $g(\frac{\partial}{\partial t},\bullet)=0$,
\item a $(1,1)$-tensor $f$ on $M\times [0,1]$ which induces on each $M_t$ an $f$-structure $f_t$ such that $f_t(\frac{\partial}{\partial t})=0$,
\item pointwise linearly independent vector fields $\{\xi_1,...,\xi_s\}$ on $M\times [0,1]$ which are tangent to $M$ (again we denote their restriction to $M_t$ by $\xi_{kt}$), 
\end{enumerate}
such that each $M_t$ with the data induced on it is an almost $\mathcal{K}$-manifold and the structure induced on $M_0$ is precisely the initial structure on $M^{2n+s}$. We say that $M_t$ is a deformation of (almost) $\mathcal{S}$-structures or (almost) $\mathcal{C}$-structures if each $M_t$ is a (almost) $\mathcal{S}$-manifold or (almost) $\mathcal{C}$-manifold respectively.
\end{defi}
\begin{rem}\begin{enumerate}
\item We note that it is sufficient to specify the data given in the above definition to define an almost $\mathcal{K}$-structure on each $M_t$ as other data (such as the $2$-form $F$) can be computed from it.
\item It is also important to note that in particular such a deformation is also a deformation of the transversely K\"ahler foliation.
\end{enumerate}
\end{rem}
We start by noting the following simple consequence of our main results:
\begin{tw}\label{bettiStab} Let $\{M^{2n+s}_t\}_{t\in [0,1]}$ be a deformation of compact almost $\mathcal{K}$-manifolds and let $b^k_t:[0,1]\to\mathbb{N}$ be the function assigning to each $t$ the $k$-th basic Betti number of $M_t$. Then:
\begin{enumerate}
\item if $M_t$ is a deformation of almost $\mathcal{S}$-manifolds such that each $M_t$ satisfies the basic hard Lefschetz property then the function $b^k_t$ is constant. In particular this is true for deformations of $\mathcal{S}$-manifolds,
\item if $M_t$ is a deformation of almost $\mathcal{C}$-manifolds then the function $b^k_t$ is constant,
\end{enumerate}
\begin{proof} Immediate since for all $t_1,t_2\in [0,1]$ the almost $\mathcal{K}$-manifolds $M_{t_1}$ and $M_{t_2}$ satisfy the assumption of Corollary \ref{TopInvS} in the first case and Corollary \ref{TopInvC} in the second case.
\end{proof}
\end{tw}
We now study the behaviour of basic Dolbeault cohomology using an approach similar to that of \cite{Noz,My2}.  We start by recalling a result from \cite{My2} which reduces the problem to proving that the spaces of complex-valued basic harmonic $k$-forms $\mathcal{H}^k_t$ of $(\mathcal{M}_t,\mathcal{F}_t)$ form a bundle over $[0,1]$.
\begin{tw}\label{IfBundle} Let $\{(M_t,\mathcal{F}_t)\}_{t\in [0,1]}$ be a smooth family of homologically orientable transversely K\"{a}hler foliations on compact manifolds such that $\mathcal{H}^k_t$ forms a smooth family of constant dimension for any $k\in\mathbb{N}$. For a fixed pair of integers $(p,q)$ the function associating to each point $s\in [0,1]$ the basic Hodge number $h^{p,q}_t$ of $(M_t,\mathcal{F}_t)$ is constant.
\end{tw}
Using this result the study can now be concluded analogously as in \cite{My2}. While the differences in the argument are scarce we present it in full for the readers convienience following closely the exposition in \cite{My2}.
\newline\indent The first step is to consider transverse $k$-forms. We denote the space of such forms by $\Omega^{T,k}$. On such forms it is natural to consider the operator $d_T:=\pi(d)$ where $\pi$ is the projection onto transverse forms given by the Riemannian metric. Its adjoint $\delta_T$ is given by the formula:
$$\delta_T:=(-1)^k*_b^{-1} d_T *_{b},$$
Which due to homological orientability coincides on basic forms with the basic coderivative $\delta_b$. This allows us to define the transverse Laplace operator in a fashion similar to \cite{E2,Noz}:
$$\Delta^T:=\sum\limits_{k=1}^s\mathcal{L}_{\xi_k}\mathcal{L}_{\xi_k}-\delta_Td_T-d_T\delta_T,$$
and similarly as in \cite{Noz,My2} we can prove the following lemma:
\begin{lem}
The operator $\Delta^T:\Omega^{k,T}\rightarrow \Omega^{k,T}$ is strongly elliptic and self-adjoint.
\begin{proof}
Around any point $x_0$ take a local coordinate chart $(z_1,...,z_s,x_1,y_1,...,x_n,y_n)$ where $\xi_k=\frac{\partial}{\partial z_k}$ and $(x_1,y_1,...,x_n,y_n)$ are transverse holomorphic coordinates such that $(\frac{\partial}{\partial x_1},\frac{\partial}{\partial y_1},...\frac{\partial}{\partial x_n},\frac{\partial}{\partial y_n})$ are orthonormal over $x_0$ and $\eta_{k}=dz_{k}+\beta_k$ for some basic forms $\beta_k$ vanishing over $x_0$. In such coordinates the principal symbol $\sigma(\delta_{T}d_T+d_T\delta_T)$ coincide with that of the Laplacian $\Delta_b$ on the planes $z_1=...=z_s=0$ (to see this note that in these coordinates $\pi(dz_k)=-\beta_k$ and so after writing the operator in local coordinates we see that aside from the parts present in $\Delta_b$ the additional components are either of degree less then $2$ or are a multiple of some $\beta_k$ (which vanish over $x_0$) and hence in either case do not contribute to the symbol over $x_0$). For $\alpha:=\sum\limits_{i=1}^s\gamma_idz_i+\sum\limits_{i=1}^n \alpha_{2i-1}dx_i+\alpha_{2i}dy_i\in T^*_{x_0}M$ let $\sigma_{\alpha}(\Delta^{T})$ be the symbol of $\Delta^T$ at $\alpha$. The symbol $\sigma_{\alpha}(\frac{\partial^2}{\partial^2 z_k})=\gamma_k^2Id_{(\Omega^{k,T})_{x_0}}$, while the symbol of $\Delta_b$ is given by $\sigma(\Delta_b)=-(\sum\limits_{i=1}^{2n}\alpha_i^2)Id_{(\Omega^{k,T})_{x_0}}$ (see \cite{V} Lemma 5.18). This shows that the symbol $\sigma_{\alpha}(\Delta^T)=||\alpha||^2Id_{(\Omega^{k,T})_{x_0}}$ and so the operator is in fact strongly elliptic.
\newline\indent Since $\delta_{T}d_T+d_T\delta_T$ is self-adjoint it suffices to prove that each $\mathcal{L}_{\xi_k}$ is skew-symmetric. For $\alpha_1,\alpha_2\in\Omega^{k,T}$ we have:
$$\mathcal{L}_{\xi_k}(\eta_1\wedge...\wedge\eta_s\wedge\alpha_1\wedge*_b\overline{\alpha_2})=\eta_1\wedge...\wedge\eta_s\wedge\mathcal{L}_{\xi_k}(\alpha_1)\wedge*_b\overline{\alpha_2}+\eta_1\wedge...\wedge\eta_s\wedge\alpha_1\wedge*_b\mathcal{L}_{\xi_k}\overline{\alpha_2},$$
since $\mathcal{L}_{\xi_k}\eta_l=0$ and $\mathcal{L}_{\xi_k}*_b=*_b\mathcal{L}_{\xi_k}$. Hence, we only need to prove that the left hand side integrates to zero over $M$. But we can write it as:
$$di_{\xi_k}(\eta_1\wedge...\wedge\eta_s\wedge\alpha_1\wedge*_b\overline\alpha_2)=d(\eta_1\wedge...\wedge\eta_{k-1}\wedge\eta_{k+1}\wedge...\wedge\eta_s\wedge\alpha_1\wedge*_b\overline\alpha_2),$$
now it suffices to note that the right hand-side is exact and hence integrates to zero.
\end{proof}
\end{lem}

With this we a ready to prove the following result:
\begin{prop}\label{BundleTrue} Let $\{M_t\}_{t\in [0,1]}$ be a smooth family of $\mathcal{C}$-manifolds or $\mathcal{S}$-manifolds over an interval. Then the spaces $\mathcal{H}^k_t$ of complex-valued basic harmonic $k$-forms on $M_t$ constitute a bundle over $[0,1]$.
\begin{proof} We start by using the results of \cite{KS} in a fashion similar to \cite{Noz} and \cite{My2} in order to contain our problem in some smooth vector bundle (with fibers of finite dimension). Using the Spectral Theorem for smooth families of strongly elliptic self-ajoint operators (see Theorem $1$ of \cite{KS}) for the family $\Delta^{k,T}_t$ we get a complete system of eigensections $\{e_{th}\}_{h\in{\mathbb{N}}, t\in[0,1]}$ together with the corresponding eigenvalues $\lambda_h(t)$ which form an ascending sequence in $[0,\infty)$ with a single accumulation point at infinity. Fix a point $t_0\in[0,1]$ and let $k_0$ be the largest number such that for $h\in\{1,...,k_0\}$ we have $\lambda_h(t_0)=0$. Consider the family of vector spaces $\mathcal{E}_t=span\{e_{th}\text{ }|\text{ }h\in\{1,...,k_0\}\}$. Since the only accumulation point of the sequence $\lambda_{h}(t_0)$ is infinity we can find a small disc around $0$ in $\mathbb{C}$ such that the only eigenvalue of $\Delta^{k,T}_{t_0}$ contained in this disc is zero. Using Theorem $2$ of \cite{KS} we establish that for each $h$ the eigenvalues $\lambda_h(t)$ form a continuous function and hence in a small neighbourhood $U$ of $t_0$ all $t\in U$ are contained in this disc as well. This allows us to conclude by using Theorem 3 of \cite{KS} that $P_{\mathcal{E}_t}(\tilde{e}_{th})$ for $h\in\{1,...,k_0\}$ form smooth sections of $\Omega^{k,T}$ over a small neighbourhood $U'\subset U$ of $t_0$ which span $\mathcal{E}_t$ (where $P_{\mathcal{E}_t}$ is the projection onto $\mathcal{E}_t$ and $\tilde{e}_{th}$ are the extensions of $e_{t_0h}$ with the use of some partition of unity over $[0,1]$). Shrinking the neighbourhood is necessary to retain linear independence of $\tilde{e}_{th}$. Hence, we have shown that $\mathcal{E}_t$ form a bundle over $U'$.
\newline\indent Now we consider the operator $\mathcal{L}_t=(\mathcal{L}_{\xi_1t},...,\mathcal{L}_{\xi_st}):\mathcal{E}_t\rightarrow \bigoplus\limits_{i=1}^s\Omega^{k,T}_t$. Note that $Ker\mathcal{L}_{t_0}=\mathcal{H}^k_{t_0}$. Via a standard rank argument there is a small neighbourhood $U''\subset U'$ of $t_0$ such that $dim(Ker\mathcal{L}_{t_0})\geq dim(Ker\mathcal{L}_{t})$. However, $Ker\mathcal{L}_{t}\supset \mathcal{H}^k_t$ and since $dim(\mathcal{H}^k_t)=dim(\mathcal{H}^k_{t_0})$ (by Theorem \ref{bettiStab}) we have the following:
$$dim(Ker\mathcal{L}_{t_0})\geq dim(Ker\mathcal{L}_{t})\geq dim(\mathcal{H}^k_t)=dim(\mathcal{H}^k_{t_0})=dim(Ker\mathcal{L}_{t_0}).$$
Hence, all of the dimensions above are equal and $Ker\mathcal{L}_{t}= \mathcal{H}^k_t$. But this implies that $\mathcal{H}^k_s$ can be described as a kernel of a morphism of bundles and since its dimension is constant we conclude that it is a bundle (over $U''$). It immediately follows that $\mathcal{H}^k_t$ forms a bundle over $[0,1]$ since it is a family of subspaces of a bundle with local trivializations around any point.
\end{proof}
\end{prop}
Combining Theorem \ref{IfBundle} and Proposition \ref{BundleTrue} we get the main result of this section:

\begin{tw}\label{hodgeStab} Let $\{M^{2n+s}_t\}_{t\in [0,1]}$ be a deformation of compact $\mathcal{C}$-manifolds or $\mathcal{S}$-manifolds. Then the function $h^{p,q}_t:[0,1]\to\mathbb{N}$, assigning to each $t$ the $(p,q)$-th basic Hodge number of $M_t$, is constant.
\end{tw}
\bigskip
\Small{\textbf{Statements and Declarations} The authors declare no competing interests.}
\bigskip

\Small{\textbf{Availability of data and material} Data sharing not applicable to this article as no datasets were generated or analysed during the current study.}

\end{document}